\documentclass[11pt, oneside]{amsart} 

\usepackage[english]{babel} 
\usepackage{graphics,color,pgf}
\usepackage{epsfig}
\usepackage[ansinew]{inputenc}
\usepackage[all]{xy}
\usepackage{tikz-cd}
\usepackage{bbm}
\usetikzlibrary{matrix,arrows,decorations.pathmorphing}
\newdir{ >}{!/8pt/@{}*@{>}}
\usepackage{amssymb, amsmath,amsthm, mathtools, amscd}
\usepackage{mathrsfs}
\usepackage[margin=1.5in]{geometry}
\usepackage{tikz}
\usetikzlibrary{automata, arrows.meta, positioning}
\usepackage[pagebackref, pdfborder={0 0 0}
  ,urlcolor=black,hypertexnames=false]{hyperref}

\usepackage[margin=1.5in]{geometry}
  
\makeatletter
\newtheorem*{rep@theorem}{\rep@title}
\newcommand{\newreptheorem}[2]{%
\newenvironment{rep#1}[1]{%
 \def\rep@title{#2 \ref{##1}}%
 \begin{rep@theorem}}%
 {\end{rep@theorem}}}
\makeatother

\newtheorem{intro_thm}{Theorem}

\theoremstyle{definition}

\theoremstyle{remark}

\newtheoremstyle{TheoremNum}
        {0.2 cm}{0.2 cm}              
        {\itshape}                      
        {}                              
        {}                     
        {.}                             
        { }                             
        {\thmname{\bfseries #1}\thmnote{ \bfseries #3}}
    \theoremstyle{TheoremNum}

\newtheorem{claim*}[rec_thm]{Claim}

\begin{document}

\title[Uniformly bounded representations into finite Von Neumann algebras]{Uniformly bounded representations of discrete measured groupoids into finite Von Neumann algebras}

\author[A. Savini]{A. Savini}
\address{Department of Mathematics, University of Milano-Bicocca, Milan, Italy}
\email{alessio.savini@unimib.it}

\date{\today.\ \copyright{\ A. Savini}.}

\maketitle

\begin{abstract}
Let $(\mathcal{G},\nu)$ be a $t$-discrete ergodic groupoid. Consider a finite Von Neumann algebra $\mathcal{M}$ with separable predual. We prove that every uniformly bounded measurable representation $\rho:\mathcal{G} \rightarrow \mathrm{GL}(\mathcal{M})$ into the invertible elements of $\mathcal{M}$ is similar to a unitary representation. This generalizes a previous result by Bates and Giordano.
\end{abstract}

\section{Introduction}

Vasilescu and Zsid\`o \cite{vasilescu} proved that every uniformly bounded subgroup of a Von Neumann algebra is actually conjugate to a unitary subgroup. This result was later generalized by Zimmer \cite{zimmer:unitary} and Bates and Giordano \cite{bates:giordano} in the context of measurable cocycles. More precisely, the latter authors proved that every uniformly bounded Borel cocycle into the invertible elements of a finite Von Neumann algebra is actually cohomologous to a unitary cocycle.

By viewing a measurable cocycle as a representation of the semidirect groupoid associated to the ergodic action of a countable group, it is natural to ask if the result by Bates and Giordano can be actually extended to the more general context of ergodic groupoids. Let $(\mathcal{G},\nu)$ be a $t$-discrete ergodic groupoid. We fix a Von Neumann algebra $\mathcal{M}$ and we denote by $\mathrm{GL}(\mathcal{M})$ the set of invertible elements in $\mathcal{M}$ with the strong operator topology. A map $\rho:\mathcal{G} \rightarrow \mathrm{GL}(\mathcal{M})$ is a \emph{morphism} if there exists a subset of units $U \subset \mathcal{G}^{(0)}$ of full measure such that the restriction to an inessential contraction $\rho:\mathcal{G}|_U^{(1)} \rightarrow \mathrm{GL}(\mathcal{M})$ is a measurable homomorphism. We say that $\rho$ is \emph{uniformly bounded} if there exists $C>0$ such that
$$
\lVert \rho(g) \rVert < C
$$
for every $g$ in some inessential contraction of $\mathcal{G}$. Two morphisms $\rho_1,\rho_2:\mathcal{G} \rightarrow \mathrm{GL}(\mathcal{M})$ are \emph{similar} if there exists a full measure subset $U \subset \mathcal{G}^{(0)}$ and a measurable map $h:\mathcal{G}^{(0)} \rightarrow \mathrm{GL}(\mathcal{M})$ such that
$$
\rho_2(g)=h(t(g))\rho_1(g)h(s(g))^{-1}
$$
makes sense and it holds for every $g \in \mathcal{G}|_U^{(1)}$. 

We prove the following:

\begin{intro_thm}\label{thm:bounded:rep}
Let $(\mathcal{G},\nu)$ be a $t$-discrete ergodic groupoid and let $\mathcal{M}$ be a finite Von Neumann algebra with separable predual. Every uniformly bounded measurable representation $\rho:\mathcal{G} \rightarrow \mathrm{GL}(\mathcal{M})$ is similar to a unitary representation.
\end{intro_thm}

We follow the same strategy of Boutonnet and Roydor \cite{boutonnet}: we exploit the geometry of the circumcenter for $\mathrm{CAT}(0)$-spaces to build the explicit similarity. 

\section{Preliminary definitions}

\subsection{Discrete ergodic groupoids} A \emph{groupoid} $\mathcal{G}$ is a small category where every morphism is invertible. The groupoid $\mathcal{G}$ is determined by both the set of objects, namely the \emph{unit space} $\mathcal{G}^{(0)}$, and by the set of morphisms, denoted by $\mathcal{G}^{(1)}$. The set of \emph{composable elements} is given by $\mathcal{G}^{(2)} \subset \mathcal{G}^{(1)} \times \mathcal{G}^{(1)}$. We denote the natural \emph{source} and \emph{target} maps by $s:\mathcal{G}^{(1)} \rightarrow \mathcal{G}^{(0)}$ and $t:\mathcal{G}^{(1)} \rightarrow \mathcal{G}^{(0)}$, respectively. Given units $x,y \in \mathcal{G}^{(0)}$, we set
$$
\mathcal{G}_x:=s^{-1}(x), \ \ \  \mathcal{G}^y:=t^{-1}(y). 
$$
By viewing $\mathcal{G}^{(0)}$ as a subset of $\mathcal{G}^{(1)}$ via the identification $x \mapsto 1_x$, one gets that $s(g)=g^{-1}g$ and $t(g)=gg^{-1}$. 
A subset $U \subset \mathcal{G}^{(0)}$ of units is called $\mathcal{G}$-\emph{invariant}, or simply \emph{invariant}, if we have that, for every $g \in \mathcal{G}^{(1)}$, $s(g) \in U$ if and only if $t(g) \in U$. Given a subset $U \subset \mathcal{G}^{(0)}$, the \emph{restricted groupoid} $\mathcal{G}|_U$ is the subgroupoid with objects $U$ and morphisms given by
$$
\mathcal{G}|_U^{(1)}:=\{ g \in \mathcal{G}^{(1)} \ | s(g),t(g) \in U\}. 
$$

A \emph{standard Borel space} $(X,\mathcal{B})$ is isomorphic to a Borel subset of a Polish space. A standard Borel space with a $\sigma$-finite Borel measure is called \emph{Lebesgue space}. It is a \emph{standard probability space} if the measure is finite and normalized to one. A groupoid $(\mathcal{G},\mathcal{B})$, or simply $\mathcal{G}$, is a \emph{Borel groupoid} if $\mathcal{G}^{(0)}$ and $\mathcal{G}^{(1)}$ are standard Borel spaces and the multiplication map $\mathcal{G}^{(2)} \rightarrow \mathcal{G}^{(1)}$ and the inverse map $\mathcal{G}^{(1)} \rightarrow \mathcal{G}^{(1)}$ are Borel. 
A \emph{measurable morphism} between two Borel groupoids $\mathcal{G}$ and $\mathcal{H}$ is a measurable map $\varphi:\mathcal{G} \rightarrow \mathcal{H}$ which is a functor between  categories. In this context, both the maps $\varphi^{(0)}:\mathcal{G}^{(0)} \rightarrow \mathcal{H}^{(0)}$ and $\varphi^{(1)}:\mathcal{G}^{(1)}\rightarrow \mathcal{H}^{(1)}$ are Borel. 

A groupoid is \emph{t}-\emph{discrete} if the target map has countable fibers. For a Borel groupoid which is $t$-discrete, given a probability measure $\mu$ on $\mathcal{G}^{(0)}$, we can define a measure $\nu$ on $\mathcal{G}^{(1)}$ as follows: for $E \subset \mathcal{G}^{(1)}$, we set
$$
\nu(E)=\int_{\mathcal{G}^{(0)}} |\mathcal{G}^x \cap E|d\mu(x). 
$$
We say that the measure $\mu$ is \emph{invariant} (or \emph{quasi-invariant}), if the measure $\nu$ (or its measure class, respectively) is invariant under the inverse map $g \mapsto g^{-1}$. When $\mu$ is quasi-invariant, we call the pair $(\mathcal{G},\nu)$ a \emph{$t$-discrete measured groupoid}. Sometimes we may need to use the triple $(\mathcal{G},\nu,\mu)$ to specify the measure defined on the units $\mathcal{G}^{(0)}$. We say that a $t$-discrete measured groupoid $(\mathcal{G},\nu,\mu)$ is \emph{ergodic}, if for every invariant subset $U \subset \mathcal{G}^{(0)}$, we have either that $\mu(U)=0$ or $\mu(\mathcal{G}^{(0)} \setminus U)=0$. 

For a measured groupoid $(\mathcal{G},\nu)$, given a full measure subset $U \subset \mathcal{G}^{(0)}$, we call \emph{inessential contraction} to $U$ the restricted groupoid $\mathcal{G}|_U$ endowed with the restricted measure. A \emph{morphism} $\varphi:\mathcal{G} \rightarrow \mathcal{H}$ from a $t$-discrete ergodic groupoid $(\mathcal{G},\nu)$ to another Borel groupoid $\mathcal{H}$ is simply a measurable morphism $\varphi:\mathcal{G}|_U \rightarrow \mathcal{H}$ defined on some inessential contraction of $\mathcal{G}$. Two morphisms $\varphi_1,\varphi_2:\mathcal{G} \rightarrow \mathcal{H}$ are \emph{similar} if there exists a measurable map $\psi:\mathcal{G}^{(0)} \rightarrow \mathcal{H}$ and there is a full measure subset $U \subset \mathcal{G}^{(0)}$ such that 
$$
s(\psi(x))=\varphi_1^{(0)}(x), \ \ \ t(\psi(x))=\varphi_2^{(0)}(x),
$$
for every $x \in U$, and
$$
\varphi_2(g)=\psi(t(g))\varphi_1(g)\psi(s(g))^{-1}
$$
for every $g \in \mathcal{G}|_U^{(1)}$. 

\subsection{Circumcenter in non-positive curvature} Given a metric space $(X,d)$, we consider a non-empty bounded subset $B$. The \emph{circumradius} of $B$ is defined as
$$
r(B):=\inf_{x \in X} \sup_{y \in B} d(x,y). 
$$
A point $x \in X$ is called a \emph{circumcenter} if the closed ball with center $x$ and radius $r(B)$ contains the set $B$. 

A geodesic metric space $(X,d)$ is $\mathrm{CAT}(0)$ if it satisfies the following inequality: for every $x_1,x_2 \in X$, there exists $z \in X$ such that for every $x \in X$ it holds
\begin{equation}\label{eq:semipar}
d(x_1,x_2)^2+4d(z,x)^2 \leq 2d(x,x_1)^2+2d(x,x_2)^2. 
\end{equation}
In a complete $\mathrm{CAT}(0)$-space every bounded set admits a unique circumcenter and it lies in the closed convex hull of the set \cite[Theorem 11.26, Theorem 11.27]{abramenko}. 

Let now $\mathcal{M}$ be a finite Von Neumann algebra with trace $\tau$. Given an element $x \in \mathcal{M}$, its $L_2$-norm is given by $\lVert x \rVert_2:=\tau(x^*x)^{1/2}$. We denote by $\mathrm{GL}^+(\mathcal{M})$ the set of positive invertible elements in $\mathcal{M}$. Given $a,b \in \mathrm{GL}^+(\mathcal{M})$ we define
$$
d(a,b):=\lVert \ln(a^{-1/2}ba^{-1/2}) \rVert_2. 
$$
This $d$ is actually a metric on $\mathrm{GL}^+(\mathcal{M})$ and satisfies the following properties \cite{miglioli}:
\begin{itemize}
\item For any $g \in \mathrm{GL}(\mathcal{M})$, we have that $d(a,b)=d(g^*ag,g^*bg)$. 
\item The metric $d$ satisfies Inequality \eqref{eq:semipar}. 
\item For every $c>1$, the metric $d$ is equivalent to the one given by the norm $\lVert \ \cdot \ \rVert_2$ on the set
$$
\mathrm{GL}_c(\mathcal{M}):=\{ x \in \mathrm{GL}^+(\mathcal{M}): 1/c \leq x \leq c \}. 
$$
For every $c>1$, the space $\mathrm{GL}_c(\mathcal{M})$ is a $\mathrm{CAT}(0)$-space which is bounded, complete and separable. 
\end{itemize}

\section{Proof of the main theorem}

We start with a measurable representation $\rho:\mathcal{G} \rightarrow \mathrm{GL}(\mathcal{M})$ which is uniformly bounded. For every $x \in \mathcal{G}^{(0)}$ we define the set 
$$
B_x:=\{ \rho(g)^\ast \rho(g) \ | \ g \in \mathcal{G}_x\}.
$$
By the uniform boundedness of $\rho$, we can assume that there exist $c>1$ and a full measure subset $U \subset \mathcal{G}^{(0)}$ such that 
$$
B_x \subset \mathrm{GL}_c(\mathcal{M}), \ \ \ \rho(g)^\ast B_{t(g)} \rho(g)=B_{s(g)},
$$
for all $x \in U$ and all $g \in \mathcal{G}|_U^{(1)}$. Let $\sigma(x) \in \mathrm{GL}_c(\mathcal{M})$ the circumcenter of the set $B_x$. As noticed by Abramenko and Brown \cite[Theorem 11.27, Lemma 11.28]{abramenko} the element $\sigma(x)$ is also the unique circumcenter of $B_x$ inside $\mathrm{GL}(\mathcal{M})$. By exploiting the proprerties of the distance $d$ on $\mathrm{GL}^+(\mathcal{M})$, together with the uniqueness of the circumcenter, we obtain that
$$
\rho(g)^\ast \sigma(t(g)) \rho(g)=\sigma(s(g)),
$$
for all $g \in \mathcal{G}|_U^{(1)}$. 

We claim that the map $x \mapsto \sigma(x)^{1/2} \in \mathcal{M}$ coincides almost everywhere with a Borel map $\psi$. We follow \cite[Lemma 3.18]{delaroche05} to prove the claim. Since $\mathcal{G}$ is $t$-discrete, for all $\vartheta \in \mathrm{GL}_c(\mathcal{M})$, the map 
$$
x \in \mathcal{G}^{(0)} \mapsto r(\vartheta, B_x):=\sup_{g \in \mathcal{G}_x} d(\vartheta,\rho(g)^\ast \rho(g))
$$
is measurable, by the measurability of $\rho$. If we fix $x \in \mathcal{G}^{(0)}$, the map $\vartheta \mapsto r(\vartheta, B_x)$ is continuous. Since $\mathrm{GL}_c(\mathcal{M})$ is separable, the map 
$$
x \in \mathcal{G}^{(0)} \mapsto r(B_x):=\inf_{\vartheta \in \mathrm{GL}_c(\mathcal{M})} r(\vartheta,B_x)
$$
is actually an infimum over a countable dense subset, and hence it is measurable. 

For any natural $n \geq 1$, we define a measurable bundle $\mathcal{D}_n$ over $\mathcal{G}^{(0)}$ as follows:
$$
\mathcal{D}_n:=\{ (x,\vartheta) \in \mathcal{G}^{(0)} \times \mathrm{GL}_c(\mathcal{M}) \ | \ r(\vartheta,B_x)^2 \leq r(B_x)^2 + 1/n \}. 
$$
By \cite[Theorem A.9]{zimmer:libro}, for every $n \geq 1$, there exists a measurable section $\xi_n:\mathcal{G}^{(0)} \rightarrow \mathrm{GL}_c(\mathcal{M})$ such that $(x,\xi_n(x)) \in \mathcal{D}_n$ for every $x \in V \subset \mathcal{G}^{(0)}$, with $V$ conull. By Inequality \eqref{eq:semipar}, there must exists $z \in \mathrm{GL}_c(\mathcal{M})$ such that
$$
d(\sigma(x),\xi_n(x))^2+4d(z,\vartheta)^2 \leq 2d(\sigma(x),\vartheta)^2+2d(\xi_n(x),\vartheta)^2
$$
for all $\vartheta \in B_x$. By taking the supremum with respect to $B_x$ we obtain that
$$
d(\sigma(x),\xi_n(x))^2 + 4r(z,B_x)^2 \leq 2r(\sigma(x),B_x)^2+2r(\xi_n(x),B_x)^2 \leq 4r(B_x)^2 +2/n.
$$
The inequality $r(z,B_x) \geq r(B_x)$ implies that $\lim_n\xi_n=\sigma$ almost everywhere and hence $x \mapsto \sigma(x)^{1/2}$ coincides almost everywhere with a Borel map $\psi$. 

Now it is sufficient to observe that 
$$
g \mapsto \psi(t(g))\rho(g)\psi(s(g))^{-1}
$$
is a unitary representation similar to $\rho$. This concludes the proof. 
\bibliographystyle{amsalpha}
\bibliography{biblionote}

\end{document}